\documentclass[12pt]{article}
\usepackage[latin1]{inputenc}
\usepackage{amssymb,amsmath,amsfonts}
\usepackage{xcolor}
\newtheorem{theorem}{Theorem}
\newtheorem{proposition}[theorem]{Proposition}
\newtheorem{lemma}[theorem]{Lemma}

\newtheorem{remark}[theorem]{Remark}

\newcommand{\R}{\mathbb{R}}
\newcommand{\U}{\mathcal{U}}
\newcommand{\V}{\mathcal{V}}
\newcommand{\Les}{\mathcal{L}}
\newcommand{\Ss}{\mathcal{S}}

\newcommand{\spa}{\mbox{span\,}}

\newcommand{\po}{{\hspace*{-1ex}}{\bf .  }}
\def\lp{{\langle\!\langle}}
\def\rp{{\rangle\!\rangle}}
\def\<{{\langle}}
\def\>{{\rangle}}
\def\Sal{{\cal S}}
\def\J{{\cal J}}
\def\T{{\cal T}}
\def\a{\alpha}
\def\be{\begin{equation} }
\def\ee{\end{equation} }

\def\proof{\noindent\emph{Proof: }}
\def\qed{\ifhmode\unskip\nobreak\fi\ifmmode\ifinner
\else\hskip5 pt \fi\fi\hbox{\hskip5 pt \vrule width4 pt
height6 pt  depth1.5 pt \hskip 1pt }}
\newcommand\blfootnote[1]{
\begingroup
\renewcommand\thefootnote{}\footnote{#1}
\addtocounter{footnote}{-1}
\endgroup
}
\begin{document}

\title{The second fundamental form of the real\\ 
Kaehler submanifolds}
\author{S. Chion and M. Dajczer}
\date{}
\maketitle

\begin{abstract} 
Let $f\colon M^{2n}\to\R^{2n+p}$, $2\leq p\leq n-1$, be an 
isometric immersion of a Kaehler manifold into Euclidean space. 
Yan and Zheng conjectured in \cite{YZ} that if the codimension 
is $p\leq 11$ then,  along any connected component of an 
open dense subset of $M^{2n}$, the submanifold is as follows: 
it is either foliated by holomorphic submanifolds 
of dimension at least $2n-2p$ with tangent spaces in
the kernel of the second fundamental form  whose images are 
open subsets of affine vector subspaces, or it is embedded 
holomorphically in a Kaehler submanifold of $\R^{2n+p}$ of larger 
dimension than $2n$. This bold conjecture was proved by Dajczer 
and Gromoll just for codimension three and then by Yan and Zheng 
for codimension four. 

In this paper we prove that the second fundamental form of the 
submanifold behaves pointwise as expected in case that the 
conjecture is true. This result is a first fundamental step for
a possible classification of the non-holomorphic Kaehler submanifolds 
lying with low codimension in Euclidean space. A counterexample shows 
that our proof does not work for higher codimension, indicating that 
proposing $p=11$ in the conjecture as the largest codimension is 
appropriate.
\end{abstract}
\blfootnote{\textup{2020} \textit{Mathematics Subject Classification}:
53B25, 53C42.}
\blfootnote{\textit{Key words}: Real Kaehler submanifolds,
the index of complex relative nullity.}

An isometric immersion $f\colon M^{2n}\to\R^{2n+p}$ is called a 
\emph{real Kaehler submanifold} if $(M^{2n},J)$ is a connected 
Kaehler manifold of complex dimension $n\geq 2$ isometrically 
immersed into Euclidean space with local substantial codimension $p$. 
The latter means that the image of $f$ restricted to any open subset 
of $M^{2n}$ does not lie in a proper affine subspace of $\R^{2n+p}$. 
Moreover, when $p$ is  even, we focus in the case in which $f$ 
restricted to any open subset of $M^{2n}$ is \emph{not} holomorphic 
with respect to any complex structure of the ambient space $\R^{2n+p}$.

Since the pioneering work by Dajczer and Gromoll \cite{DG1}, there
has been an increasing interest in the study of the real Kaehler 
submanifolds. The reason, in good part, it is due because when 
these submanifolds are minimal then they enjoy several of the feature 
properties of minimal surfaces. For instance, they admit 
an associated one-parameter family of non congruent isometric
minimal submanifolds all with the same Gauss map. Another one 
is being the real part of 
its holomorphic representative. Moreover, the immersions are 
pluriharmonic maps and, in some cases, they admit a Weierstrass type 
representation. For a partial account of results on this subject 
of research, as well as many references, we refer to \cite{DT}.
\vspace{1ex}

There is plenty of knowledge on real Kaehler submanifolds  
$f\colon M^{2n}\to\R^{2n+p}$ when the codimension is as low as 
$p=1,2$. For instance, in the hypersurface case, there is the 
local parametric classification obtained in \cite{DG1} that 
can be seen in \cite{DT} as Theorem~$15.14$. The classification 
of the metrically complete submanifolds with codimension 
$p=2$ follows from \cite{DG2} and \cite{FZ2}. Moreover, 
for both codimensions, the submanifolds carry a foliation by 
complex relative nullity leaves of dimension $2n-2p$ as 
described next.

Let $f\colon M^{2n}\to\R^{2n+p}$ be a real Kaehler submanifold
and let $L\subset N_fM(x)$ be a normal vector subspace at 
$x\in M^{2n}$. We denote the $\a_L\colon TM\times TM\to L$
the $L$-component of its normal vector valued second fundamental 
form $\a\colon TM\times TM\to N_fM$ and by 
$\mathcal{N}(\a_L)\subset T_xM$ the tangent vector subspace 
$$
\mathcal{N}(\a_L)=\{Y\in T_xM\colon \a_L(X,Y)=0
\;\mbox{for any}\;X\in T_xM\}.
$$
Then $\Delta(x)=\mathcal{N}(\a_{N_fM(x)})$ is called the relative 
nullity subspace of $f$ at $x\in M^{2n}$. Its complex part 
$\Delta_c(x)=\Delta(x)\cap J\Delta(x)$ is named  the 
\emph{complex relative nullity} subspace whose dimension 
$\nu_f^c(x)$ is the \emph{index of complex relative nullity}.
It is well-known that the vector subspaces $\Delta_c(x)$ 
form a smooth integrable distribution on any open subset of 
$M^{2n}$ where $\nu_f^c(x)$ is constant. Moreover, the 
totally geodesic leaves are holomorphic submanifolds of 
$M^{2n}$ as well as open subsets of even dimensional affine 
vector subspaces of $\R^{2n+p}$.
\vspace{1ex}

Real Kaehler submanifolds in codimension at least three can 
be obtained just by considering holomorphic submanifolds of 
a given real Kaehler submanifold. More precisely, let 
$F\colon N^{2n+2m}\to\R^{2n+p}$, $m\geq 1$, be a real Kaehler 
submanifold, and then let $j\colon M^{2n}\to N^{2n+2m}$ be any 
holomorphic isometric immersion. Then the \emph{composition} 
of isometric immersions $f\colon M^{2n}\to\R^{2n+p}$ given by
\be\label{composition}
f=F\circ j\colon M^{2n}\to\R^{2n+p}
\ee 
is a real Kaehler submanifold. 

It is clearly relevant to establish conditions asserting that  
a real Kaehler submanifold is locally a composition as in 
\eqref{composition}. This was achieved for $p=3$ by Dajczer and 
Gromoll \cite{DG3} and for $p=4$ by Yan and Zheng \cite{YZ} under 
the assumption that the index of complex relative nullity of 
$f\colon M^{2n}\to\R^{2n+p}$ satisfies $\nu_f^c(x)<2n-2p$ at 
any $x\in M^{2n}$.  The result in the latter paper was 
complemented by us in \cite{CDa2}.

A bold conjecture by Yan and Zheng in \cite{YZ} states, under 
the same assumption as above on the index of complex relative 
nullity, that any real Kaehler submanifold in codimension 
$p\leq 11$ is a composition as in \eqref{composition}
along connected components of an open dense subset of $M^{2n}$. 
The purpose of this paper is to walk a fundamental step in order 
to treat that rather challenging conjecture. We prove that the 
second fundamental form of the submanifold behaves pointwise as 
expected if the conjecture were true. Moreover, we have that our 
proof fails for $p=12$, indicating that proposing $p=11$ in 
the conjecture as the largest codimension seems appropriate. 
For codimension  $p\leq 6$ our result was obtained in \cite{CG} 
up to some inconsistencies in the argument; see Remark \ref{ag}. 
As for higher codimension, it is shown by this paper that the 
proof is much more difficult.
\vspace{1ex}

Before stating our main theorem we roughly explain why this result 
turns out to be the one we expected.  For this purpose, let 
$f=F\circ j\colon M^{2n}\to\R^{2n+p}$ be a composition, as in 
\eqref{composition} but where $F$ itself  is not such a composition. 
Then the second fundamental form $\a^f\colon TM\times TM\to N_f M$ 
of $f$ splits as the sum of the second fundamental forms $\a^j$ 
of $j$ and $\a^F$ of $F$ that one restricted to $TM$, and where both 
components need to satisfy certain conditions now discussed. On one hand, 
there is a vector bundle isometry $\J\in\Gamma(Aut(\Omega))$ such that
\be\label{maincondition}
\J\a^f_\Omega(X,Y)=\a^f_\Omega(X,JY)\;\;\mbox{for any}
\;\;X,Y\in\mathfrak{X}(M)
\ee 
where $\Omega=F_*N_jM$.
In fact, being $j$ holomorphic we have that $j_*JX=J^Nj_*X$ for 
any  $X\in\mathfrak{X}(M)$. Differentiating once and then taking 
normal component yields $J^N\a^j(X,Y)=\a^j(X,JY)$.  
Then $\J F_*|_{N_jM}=F_*J^N|_{F_*N_jM}$ satisfies the requirement. 
On the other hand, since $F$ is not a composition then $\a^F$ should 
have a large index of complex relative nullity, and hence the same  
remains to be the case when it is restricted to~$TM$. 
\vspace{1ex}

Let $N_1(x)\subset N_fM(x)$ denote the vector subspace 
spanned at $x\in M^{2n}$ by the second fundamental form
of $f$, namely, $N_1(x)=\spa\{\a(X,Y)\colon X,Y\in T_xM\}$. 
It is usually called the first normal space of $f$ at 
$x\in M^{2n}$. 
Then let $Q(x)\subset N_1(x)$ the complex vector subspace 
defined as
$$
Q(x)=\{\eta\in N_1(x)\colon\<\eta,\a(Z,T)\>
=\<\bar\eta,\a(Z,JT)\>\;\mbox{for any}\; Z,T\in T_xM\}
$$
where if $\eta=\sum_{i=1}^k\a(X_i,Y_i)$ then $\bar\eta
=\sum_{i=1}^k\a(X_i,JY_i)$. 

\begin{theorem}\po\label{main} Let $f\colon M^{2n}\to\R^{2n+p}$, 
$2\leq p\leq n-1$, be a real Kaehler submanifold whose index of 
complex relative nullity satisfies $\nu_f^c(x_0)<2n-2p$ at a
point $x_0\in M^{2n}$. If $p\leq 11$ then the 
following facts hold:
\begin{itemize}
\item[(i)] If $Q=Q(x_0)$ then $\dim Q=\ell>0$ and there is 
an isometry $\J\in Aut(Q)$ such that
$$
\J\a_Q(X,Y)=\a_Q(X,JY)\;\;\mbox{for any}\;\; X,Y\in T_{x_0}M.
$$  
\item[(ii)] If $N_1(x_0)=Q\oplus P$ is an orthogonal decomposition 
then $\nu^c(\a_P)\geq 2(n-p+\ell)$.
\end{itemize}
\end{theorem}

If the submanifold satisfies $\dim N_1(x_0)=q<p$ then the 
proof of Theorem \ref{main} gives a stronger result.
Indeed, one can replace  the assumption $p\leq 11$ by  
$q\leq 11$ and assume in part $(ii)$ that $\nu_f^c(x_0)<2n-2q$.
\vspace{1ex}

The extrinsic assumption on index of complex relative nullity in 
Theorem \ref{main} can be replaced by an intrinsic hypothesis,
namely, there is no complex vector subspace $L^{2n-2p}\subset T_{x_0}M$ 
such that the sectional curvature satisfies $K_M(P)=0$ for any
plane $P^2\subset L^{2n-2p}$.  Notice that part $(i)$ gives that 
$\J$ is a complex structure, that is, that we have $\J^2=-I$. 
It also yields that $\a_Q(JX,Y)=\a_Q(X,JY)$ holds for any 
$X,Y\in T_{x_0}M$.  Finally, we observe that the inequality 
$\nu_f^c(x)<2n-2p$ holds in a neighborhood of $x_0$ in $M^{2n}$.
\vspace{1ex}

Although the above result can be seen as 
a validation of the Yan and Zheng conjecture at the level of the structure 
of the second fundamental form of the submanifold, it is a distance apart 
from proving that the conjecture is true. In fact, we believe that for 
codimensions $p\geq 7$ there is just one other possibility, namely, 
that we may have complex ruled submanifolds that are not compositions. 
By being complex ruled we mean that there is a holomorphic foliation of 
$M^{2n}$ such that the image by $f$ of each leaf is part of an affine 
vector subspace of $\R^{2n+p}$ but it does not have to be part of the 
complex relative nullity.
\vspace{1ex}

An immediate application of Theorem \ref{main} is the following 
result under a pinching curvature condition.

\begin{theorem}\po\label{cor} Let $f\colon M^{2n}\to\R^{2n+p}$, 
$2\leq p\leq n-1$, be a real Kaehler submanifold whose index of 
complex relative nullity satisfies $\nu_f^c(x_0)<2n-2p$ at 
$x_0\in M^{2n}$. If $p\leq 11$ there is a neighborhood $U$ of 
$x_0$ such that an any point $x\in U$ there is a complex vector 
subspace $L^{2m}\subset T_xM$ with $m\geq n-p+\ell$ where 
$\dim Q(x)=\ell>0$ such that for any complex plane 
$P^2\subset L^{2m}$ the sectional curvature satisfies 
$K_M(P)\leq 0$.
\end{theorem}

For $p\leq n$ and without the assumption on the index of 
complex relative nullity the weaker estimate $m\geq n-p$ 
was given as Corollary $15.6$ in \cite{DT}.
\vspace{1ex}

Finally, we observe that if $p$ is even and  
$f\colon M^{2n}\to\R^{2n+p}$ is holomorphic with respect to 
some complex structure in the ambient space then $Q(x)=N_1(x)$ 
holds everywhere and hence both results given above are trivial.

\section{Preliminaries}

This section provides several basic facts used throughout 
the paper. 
\vspace{2ex}

Let $\varphi\colon V_1\times V_2\to W$ 
denote a bilinear form between real vector spaces of finite dimension. 
The image of $\varphi$ is the vector subspace of $W$ defined by
$$
\mathcal{S}(\varphi)
=\spa\{\varphi(X,Y)\;\mbox{for all}\; X\in V_1\;\mbox{and}\;Y\in V_2\}
$$
whereas the (right) \emph{nullity} of $\varphi$ is the vector 
subspace of $V_2$ given by
$$
\mathcal{N}(\varphi)=\{Y\in V_2\colon\varphi(X,Y)=0
\;\mbox{for all}\;X\in V_1\}
$$
whose dimension $\nu(\varphi)$ is the index of nullity of $\varphi$.
\vspace{1ex}

A vector $X\in V_1$ is called a (left) regular element of 
$\varphi$ if $\dim\varphi_X(V_2)=\kappa(\varphi)$ where
$$
\kappa(\varphi)=\max_{X\in V_1}\{\dim\varphi_X(V_2)\}
$$
and $\varphi_X\colon V_2\to W$ is the linear map defined by 
$$
\varphi_XY=\varphi(X,Y).
$$
Then $RE(\varphi)\subset V_1$ 
denotes the subset of regular elements of $\varphi$. 
Given $X\in RE(\varphi)$ then the vector subspace 
$N(X)=\ker\varphi_X$ satisfies
\be\label{general}
\dim N(X)=\dim V_2-\kappa(\varphi).
\ee

Let $W$ be endowed with an inner product of any 
signature. Then we denote
$$
\U(X)=\varphi_X(V_2)\cap\varphi_X(V_2)^\perp,
$$
$\tau_\varphi(X)=\dim\U(X)$ and 
$\tau(\varphi)=\min_{X\in RE(\varphi)}\{\tau_\varphi(X)\}$.
\vspace{1ex}

The following result will be used throughout the paper without 
further reference.

\begin{proposition}\po The following facts hold:
\begin{itemize}
\item[(i)] The subset $RE(\varphi)\subset V_1$ is open 
and dense.
\item[(ii)]
If $V_1=V_2=V$ and $\varphi$ is symmetric then 
$$ 
RE^*(\varphi)=\{X\in RE(\varphi):\varphi(X,X)\neq 0\}
$$ 
is an open dense subset of $V$.
\item[(iii)] If $W$ is endowed with an inner product then
$$
RE^\#(\varphi)=\{X\in RE(\varphi):\tau_\varphi(X)=\tau(\varphi)\}
$$ 
is an open dense subset of $V_1$. 
\end{itemize}
\end{proposition}

\proof Part $(i)$ is Proposition $4.4$ in \cite{DT} whereas 
the proof of Lemma $2.1$ in \cite{DR} gives part $(iii)$. 
An easy argument gives part $(ii)$, for instance, see the 
proof of Lemma $4.5$ in \cite{DT}.
\vspace{2ex}\qed

Let $W$ be endowed with the inner product $\lp\,,\,\rp$. Then 
the bilinear form $\varphi$ is said to be \emph{flat} if 
$$
\lp\varphi(X,Y),\varphi(Z,T)\rp-\lp\varphi(X,T),\varphi(Z,Y)\rp=0
$$
for any $X,Z\in V_1$ and $Y,T\in V_2$. It is said that $\varphi$ 
is \emph{null} if 
$$
\lp\varphi(X,Y),\varphi(Z,T)\rp=0
$$ 
for any $X,Z\in V_1$ and $Y,T\in V_2$. \vspace{1ex}

Given $X\in RE(\varphi)$ we denote
$$ 
\Les(X)=\Sal(\varphi|_{V_1\times N(X)}). 
$$
Then let $\sigma_\varphi(X)=\dim\Les(X)$  and 
$\sigma(\varphi)=\min_{X\in RE(\varphi)}\{\dim\sigma_\varphi(X)\}$.

\begin{proposition}\po\label{firstst}  If $X\in RE(\varphi)$ then
$\Les(X)\subset\varphi_X(V_2)$. Moreover, if $\varphi$ is flat then
\be\label{secondst}
\Les(X)\subset\U(X)
\ee
and thus $\sigma(\varphi)\leq\sigma_\varphi(X)\leq\tau_\varphi(X)$. 
\end{proposition}

\proof See Proposition $4.6$ in \cite{DT}.\vspace{2ex}\qed

Let $U^p$ be a $p$-dimensional vector space induced with a positive 
definite inner product $\<\,,\,\>$. Set $W^{p,p}=U^p\oplus U^p$ 
and let $\pi_1\colon W^{p,p}\to U^p$ (respectively, $\pi_2$)  
denote taking the first (respectively, 
second) component of $W^{p,p}$. Then let $W^{p,p}$ be endowed with 
the inner product $\lp\,,\,\rp$ of signature $(p,p)$ given by
$$
\lp(\xi_1,\xi_2),(\eta_1,\eta_2)\rp
=\<\xi_1,\eta_1\>-\<\xi_2,\eta_2\>.
$$
Then $\mathcal{T}\in\text{Aut}(W)$ defined by 
\be\label{mathcalT}
\mathcal{T}(\xi,\eta)=(\eta,-\xi)
\ee
is a complex structure which means that $\mathcal{T}^2=-I$. 
Moreover, it holds that
$$
\lp\mathcal{T}\delta,\nu\rp=\lp\delta,\mathcal{T}\nu\rp.
$$

A vector subspace $L\subset W^{p,p}$ is called \emph{degenerate} 
if $L\cap L^\perp\neq 0$ and  \emph{nondegenerate} if otherwise.  
A degenerate vector subspace $L\subset W^{p,p}$ is called 
\emph{isotropic} if $L=L\cap L^\perp$. 

\begin{proposition}\po\label{decomp} Given a vector subspace 
$L\subset W^{p,p}$ there is a direct sum decomposition 
\be\label{thedecomp}
W^{p,p}=\U^r\oplus\hat{\U}^r\oplus\mathcal{V}^{p-r,p-r}
\ee
where $\U^r=L\cap L^\perp$, the vector subspace $\hat{\U}^r$ is 
isotropic, the vector subspace $\U^r\oplus\hat{\U}^r$ is 
nondegenerate  and $L\subset\U^r\oplus\mathcal{V}^{p-r,p-r}$ 
where $\mathcal{V}^{p-r,p-r}=(\U^r\oplus\hat{\U}^r)^\perp$.
\end{proposition}

\proof See Sublemma $2.3$ in \cite{CD} or Corollary $4.3$ 
in \cite{DT}.\qed

\begin{remark}\po {\em In the decomposition \eqref{thedecomp} 
only $\U^r$ is completely determined by $L$. 
In fact, if $\hat{\U}^r=\spa\{\xi_1,\ldots,\xi_r\}$ then
any alternative description is as
$\spa\{\xi_1+\delta_1,\ldots,\xi_r+\delta_r\}$ where 
$\{\delta_1,\ldots,\delta_r\}$ is any set of vectors 
belonging to $\mathcal{V}^{p-r,p-r}$ that span an isotropic 
subspace.
}\end{remark}
 
Let the vector space $V_2$ carry a complex structure
$J\in\text{Aut}(V_2)$. 
It is a standard fact that $V_2$ is even-dimensional and 
admits a basis of the form $\{X_j,JX_j\}_{1\leq j\leq n}$.
Assume that the bilinear form 
$\varphi\colon V_1\times V_2\to W^{p,p}$ satisfies that
\be\label{condo}
\mathcal{T}\varphi(X,Y)=\varphi(X,JY)\;\;\mbox{for any}
\;\;X\in V_1\;\mbox{and}\;Y\in V_2
\ee
and let $W^{p,p}=\U\oplus\hat\U\oplus\mathcal{V}$ be the 
decomposition given by \eqref{thedecomp} for $L=\Sal(\varphi)$. 
Then we have $\mathcal{T}\U=\U$.
In effect, if $\lp\varphi(X,Y),(\xi,\bar\xi)\rp=0$ for any 
$X\in V_1$ and $Y\in V_2$ then 
$$
\lp\varphi(X,Y),\mathcal{T}(\xi,\bar\xi)\rp
=\lp\mathcal{T}\varphi(X,Y),(\xi,\bar\xi)\rp
=\lp\varphi(X,JY),(\xi,\bar\xi)\rp=0.
$$

\begin{proposition}\po\label{even}
The following facts hold:
\begin{itemize}
\item[(i)] 
$\mathcal{T}|_{\Sal(\varphi)}\in\text{Aut}\,\,(\Sal(\varphi))$
and $\mathcal{T}|_{\U}\in\text{Aut}(\U)$ are complex structures.
\item[(ii)] The vector subspaces $\Sal(\varphi)$ and $\U$ of 
$W^{p,p}$ have even dimension.
\item[(iii)] The vector subspace $\mathcal{N}(\varphi)\subset V_2$ 
is $J$-invariant and thus of even dimension.
\item[(iv)] If $\Omega=\pi_1(\U)$ then $\dim\U=\dim\Omega$ and 
if $\varphi_\Omega=\pi_{\Omega\times\Omega}\circ\varphi$
then $\Sal(\varphi_\Omega)=\U$.
\end{itemize}
\end{proposition}

\proof The considerations given above yield parts $(i)$ to $(iii)$. 
Being the subspace $\U$ isotropic then $\pi_1|_{\U}\colon\U\to\Omega$ 
is an isomorphism. Since $\mathcal{T}\U=\U$ gives that $\pi_2(\U)=\Omega$ 
then part $(iv)$ follows.\qed

\begin{proposition}\po Let the bilinear form 
$\varphi\colon V^{2n}\times V^{2n}\to W^{p,p}$ be symmetric and 
satisfy the condition \eqref{condo}. Then
\be\label{estpluri}
4\dim\Sal(\varphi)\leq\kappa(\varphi)(\kappa(\varphi)+2).
\ee 
\end{proposition}

\proof Since $\mathcal{T}|_{\varphi_X(V)}$ is a complex structure 
then $\kappa(\varphi)=2m$. Fix $X\in RE^*(\varphi)$ and let  
$\{X_j,JX_j\}_{1\leq j\leq n}$ be a basis of $V^{2n}$ with 
$X_1=X$ such that
$$
\varphi_X(V)
=\spa\{\varphi_XX_j,\varphi_XJX_j, 1\leq j\leq m\}
$$
and $X_r,JX_r\in\ker\varphi_X$ for $r\geq m+1$.
Since Proposition \ref{firstst} yields 
$\Sal(\varphi|_{V\times \ker\varphi_X})
\subset\varphi_X(V)$ then given $Z\in V^{2n}$ and $q\geq m+1$ 
there is $Y\in\spa\{X_j,JX_j, 1\leq j\leq m\}$ such that
$$
\varphi(Z,X_q)=\varphi(X_1,Y)\;\;
\mbox{and}\;\;\varphi(Z,JX_q)=\varphi(X_1,JY).
$$
Being $\varphi$ symmetric, we have 
$$
\varphi(X,JY)=\mathcal{T}\varphi(X,Y)=\mathcal{T}\varphi(Y,X)
=\varphi(Y,JX)=\varphi(JX,Y)
$$
for any $X,Y\in V^{2n}$. Hence
$$
\Sal(\varphi)
=\spa\{\varphi(X_i,X_j),\varphi(X_i,JX_j),1\leq i\leq j\leq m\},
$$
and \eqref{estpluri} follows.\qed

\section{The proofs}

In this section we first give a general result in the theory 
of flat bilinear forms tailored for our purposes in 
this paper. After that we prove both results that 
have been stated in the Introduction.
\vspace{1ex}

Let $\a\colon V^{2n}\times V^{2n}\to U^p$ be a symmetric 
bilinear form and $J\in\text{Aut}(V)$ a complex structure. 
Then let $\gamma\colon V^{2n}\times V^{2n}\to W^{p,p}$ be the 
associated bilinear form defined by 
\be\label{gamma}
\gamma(X,Y)=(\a(X,Y),\a(X,JY)).
\ee
Then $\gamma$ is symmetric if and only if
$\a$ is \emph{pluriharmonic} with the latter meaning that  
$$
\a(JX,Y)=\a(X,JY)\;\;\mbox{for any}\;\; X,Y\in V^{2n}.
$$
If $\mathcal{T}\in\text{Aut}(W)$ is the complex structure 
given by \eqref{mathcalT} then
\be\label{T1}
\mathcal{T}\gamma(X,Y)=\gamma(X,JY)
\;\;\mbox{for any}\;\; X,Y\in V^{2n}
\ee
and thus Proposition \ref{even} applies to $\gamma$.

Let $\beta\colon V^{2n}\times V^{2n}\to W^{p,p}$ be the 
bilinear form defined by
\begin{align}\label{beta}
\beta(X,Y)
&=\gamma(X,Y)+\gamma(JX,JY)\\
&=(\a(X,Y)+\a(JX,JY),\a(X,JY)-\a(JX,Y))\nonumber.
\end{align}
By \eqref{T1} we have that 
\be\label{T2}
\mathcal{T}\beta(X,Y)=\beta(X,JY)
\;\;\mbox{for any}\;\; X,Y\in V^{2n}
\ee
and hence Proposition \ref{even} applies to $\beta$. 
Then part $(iii)$ gives that $\nu(\beta)$ is even. We observe that 
$\nu(\beta)$ was called in \cite{FZ1} the index of pluriharmonic 
nullity since it satisfies
$$
\mathcal{N}(\beta)=\{Y\in V^{2n}:\a(X,JY)
=\a(JX,Y)\;\mbox{for all}\;X\in V^{2n}\}.
$$

\begin{theorem}\po\label{alglemma2} Let the bilinear forms 
$\gamma,\beta\colon V^{2n}\times V^{2n}\to W^{p,p}$,
$p\leq n$, be flat and satisfy
\be\label{productflat2}
\lp\beta(X,Y),\gamma(Z,T)\rp=\lp\beta(X,T),\gamma(Z,Y)\rp\;\;
\mbox{for any}\;\; X,Y,Z,T\in V^{2n}.
\ee
If $p\leq 11$ and $\nu(\gamma)<2n-\dim\Sal(\gamma)$
then  the vector subspace $\Sal(\gamma)$ is degenerate.
Moreover, if $\Omega=\pi_1(\Sal(\gamma)\cap\Sal(\gamma)^\perp)$
and $U^p=\Omega\oplus P$ is an orthogonal decomposition
then the following holds:
\begin{itemize}
\item[(i)] There is an isometric complex structure 
$\J\in\text{End}(\Omega)$ so that $\a_\Omega=\pi_\Omega\circ\a$ 
satisfies 
$$
\J\a_\Omega(X,Y)=\a_\Omega(X,JY)\;\;\mbox{for any}\;\;X,Y\in V^{2n}.
$$  
\item[(ii)] The  bilinear form $\gamma_P=\pi_{P\times P}\circ\gamma$ 
is flat, the vector subspace $\Sal(\gamma_{P})$ is nondegenerate 
and $\nu(\gamma_P)\geq 2n-\dim\Sal(\gamma_P)$.
\end{itemize}
\end{theorem}

The proof of Theorem \ref{alglemma2} will require several lemmas. 

\begin{lemma}\po\label{mainlemmabis}
Let the bilinear form $\gamma\colon V^{2n}\times V^{2n}\to W^{p,p}$ 
be symmetric and flat. If $p\leq 11$ and $\Sal(\gamma)=W^{p,p}$ then 
\be\label{estim1mlb}
\nu(\gamma)\geq 2n-\kappa(\gamma)-\sigma(\gamma)
\geq 2n-2p.
\ee
\end{lemma}

\proof We argue for the most difficult case $p=11$ being the 
other cases similar but easier as $p$ decreases. The first inequality in 
\eqref{estim1mlb} just means that 
$$
\nu(\gamma)\geq 2n-\kappa(\gamma)-\sigma_\gamma(X)
\;\;\mbox{for any}\;\; X\in RE(\gamma).
$$
Thus for what follows we fix $X\in RE(\gamma)$ and prove 
the latter. Proposition \ref{decomp} yields 
\be\label{decompmlb}
W^{p,p}
=\U^\tau(X)\oplus\hat\U^\tau(X)\oplus\mathcal{V}^{p-\tau,p-\tau}(X)
\ee
where $\U^\tau(X)=\gamma_X(V)\cap\gamma_X(V)^\perp$,
$\gamma_X(V)\subset\U^\tau(X)\oplus\mathcal{V}^{p-\tau,p-\tau}(X)$
and $\tau=\tau_\gamma(X)$ for simplicity. Thus
$\kappa(\gamma)\leq 2p-\tau$. Then  \eqref{secondst} yields
$\kappa(\gamma)+\sigma_\gamma(X)\leq \kappa(\gamma)+\tau\leq 2p$.
which gives the second inequality in \eqref{estim1mlb}.

The vector subspace $\U^\tau(X)$ is zero or is by Proposition \ref{even} 
isotropic of even dimension. It follows 
from \eqref{secondst} that $0\leq\sigma\leq\tau\leq 10$  where  
$\sigma=\sigma_\gamma(X)$ for simplicity of notation.

If $\sigma=0$, that is, we have $N(X)=\mathcal{N}(\gamma)$ 
and then \eqref{estim1mlb} is just \eqref{general}. Hence we assume 
$\sigma>0$. Moreover, using first part $(iii)$ and then part $(ii)$ 
of Proposition \ref{even} we obtain that $\sigma$ is even. Thus, 
henceforth we assume $\sigma\geq 2$.
\vspace{1ex}

In view of \eqref{secondst} there is a decomposition
\be\label{decomp2mlb}
\U^\tau(X)\oplus\hat\U^\tau(X)
=\Les^{\sigma}(X)\oplus\hat\Les^{\sigma}(X)
\oplus\mathcal{V}_0^{\tau-\sigma,\tau-\sigma}
\ee
where $\hat\Les^{\sigma}(X)\subset\hat\U^\tau(X)$ is such that 
the vector subspace $\mathcal{V}_0^{\tau-\sigma,\tau-\sigma}
=(\Les^{\sigma}(X)\oplus\hat\Les^{\sigma}(X))^\perp$ is 
nondegenerate. We denote 
$\hat\gamma=\pi_{\hat\Les^{\sigma}(X)}\circ\gamma$ and show that 
$\mathcal{T}\hat\gamma(Y,Z)=\hat\gamma(Y,JZ)$, that is, that
\be\label{partv} 
\mathcal{T}|_{\Sal(\hat{\gamma})}\hat{\gamma}_YZ=\hat{\gamma}_YJZ
\;\;\mbox{for any}\;\;Y,Z\in V^{2n}.
\ee
Hence $\mathcal{T}|_{\Sal(\hat{\gamma})}$ 
is a complex structure and $\kappa_0=\kappa(\hat\gamma)$ is even.
Part $(iii)$ of Proposition \ref{even} gives that $N(X)$
is $J$-invariant. If $(\xi,\bar\xi)\in\mathcal{L}^{\sigma}(X)$
then part $(i)$ of Proposition \ref{even}
applied to $\varphi=\gamma|_{V\times N(X)}$ yields that
$\mathcal{T}(\xi,\bar\xi)\in\mathcal{L}^{\sigma}(X)$. 
Using \eqref{decompmlb} and \eqref{decomp2mlb} we have
\begin{align*}
\lp\mathcal{T}\hat\gamma(Y,Z),(\xi,\bar\xi)\rp
&=\lp\hat\gamma(Y,Z),\mathcal{T}(\xi,\bar\xi)\rp
=\lp\gamma(Y,Z),\mathcal{T}(\xi,\bar\xi)\rp
=\lp\gamma(Y,JZ),(\xi,\bar\xi)\rp\\
&=\lp\hat\gamma(Y,JZ),(\xi,\bar\xi)\rp
\end{align*}
for any $(\xi,\bar\xi)\in\mathcal{L}^{\sigma}(X)$
and this gives \eqref{partv}.

We have that 
\be\label{sigmatau}
\sigma\leq\tau(\gamma).
\ee
In effect, it follows from \eqref{general} that the dimension 
of $N(Y)$ on $RE(\gamma)$ is constant. Then by continuity 
$\sigma\leq\sigma_\varphi(Y)$ in a neighborhood of $X$ in 
$RE(\gamma)$. On the other hand, we obtain from \eqref{secondst} 
that $\sigma_\varphi(Y)\leq\tau(\gamma)$ for any 
$Y\in RE^\#(\gamma)$ which is open and dense in $V^{2n}$. 
Then \eqref{sigmatau} follows. 
\vspace{1ex}

\noindent{\it Claim.} 
Given $Z\in V^{2n}$ then $\dim\gamma_Z(N(X))$ is even and
\be\label{theclaim}
\dim\gamma_Z(N(X))\leq p-\kappa_0-\tau(\gamma)+\sigma\leq p-\kappa_0.
\ee
That $\dim\gamma_Z(N(X))$ is even follows from parts $(ii)$ and 
$(iii)$ of Proposition \ref{even} whereas \eqref{sigmatau} yields 
the second inequality in \eqref{theclaim}.

To prove the first inequality in \eqref{theclaim} it suffices to 
argue for $Z\in RE^\#(\gamma)\cap RE(\hat\gamma)$ since this subset 
of $V^{2n}$ is open and dense. Let $V_0\subset V^{2n}$ be the vector 
subspace $V_0=\gamma_Z^{-1}(\mathcal{L}^{\sigma}(X))$ and 
$s_0=\dim\gamma_Z(V_0)$. Since $N(X)\subset V_0$ by \eqref{secondst} 
then $r\leq s_0$ where $r=\dim\gamma_Z(N(X))$. Because
$Z\in RE(\hat\gamma)$ there is a vector subspace 
$V_1^{\kappa_0}\subset V^{2n}$ satisfying 
$\hat\gamma_Z(V_1)=\hat\gamma_Z(V)$. Since any vector in 
$\gamma_Z(V_1)$ has a nonzero  $\hat\Les^{\sigma}(X)$-component 
then 
\be\label{disjoint}
\gamma_Z(V_0)\cap\gamma_Z(V_1)=0.
\ee

Let $Y_0\in V_0$ satisfy $\gamma_ZY_0\in\U^{\bar\tau}(Z)$
where $\bar\tau=\tau(\gamma)$ for simplicity of notation. 
Since $\gamma_Z(V_0)\subset\Les^{\sigma}(X)$ and 
$\hat\gamma_Z(V_1)\subset \hat\Les^{\sigma}(X)$
then using \eqref{decomp2mlb} we have
$$
\lp\gamma_ZY_0,\hat\gamma_Z(V_1)\rp
=\lp\gamma_ZY_0,\gamma_Z(V_1)\rp=0.
$$
Hence 
\be\label{firstu}
\dim\gamma_Z(V_0)\cap\U^{\bar\tau}(Z)\leq\sigma-\kappa_0.
\ee

Let $Y_1\in V_1^{\kappa_0}$ satisfy $\gamma_ZY_1\in\U^{\bar\tau}(Z)$.
Since $\gamma_Z(V_0)\subset\Les^{\sigma}(X)$ and 
$\hat\gamma_Z(V_1)\subset\hat\Les^{\sigma}(X)$ 
then \eqref{decomp2mlb} gives
$$
\lp\hat\gamma_ZY_1,\gamma_Z(V_0)\rp
=\lp\gamma_ZY_1,\gamma_Z(V_0)\rp=0.
$$
Hence $\dim\pi_{\hat{\Les}(X)}(\gamma_Z(V_1)\cap\U^{\tilde{\tau}}(Z))
\leq\sigma-s_0$. But since $V_1^{\kappa_0}$ has been chosen to
satisfy that $\pi_{\hat{\Les}(X)}|_{\gamma_Z(V_1)}$ is injective, 
then 
\be\label{secondu}
\dim\gamma_Z(V_1)\cap\U^{\bar\tau}(Z)\leq\sigma-s_0.
\ee

The decomposition \eqref{decompmlb} for $Z$ yields 
$\gamma_Z(V)\subset\U^{\bar\tau}(Z)\oplus 
\mathcal{V}^{p-\bar\tau,p-\bar\tau}(Z)$. Let the vector 
subspace $\mathcal{R}\subset\gamma_Z(V)$ be such that
$\gamma_Z(V)=(\gamma_Z(V)\cap\U^{\bar{\tau}}(Z))\oplus\mathcal{R}$.
Since any vector in $\mathcal{R}$ has a nonzero 
$\mathcal{V}^{p-\bar\tau,p-\bar\tau}(Z)$-component
then $\pi_{\mathcal{V}(Z)}|_{\mathcal{R}}$ is injective. 

Set $\Ss=\pi_{\mathcal{V}(Z)}(\gamma_Z(V_0)\cap\mathcal{R})$ and 
$\hat\Ss=\pi_{\mathcal{V}(Z)}(\gamma_Z(V_1)\cap\mathcal{R})$. Since 
$\dim\gamma_Z(V_0)=s_0$ and $\dim\gamma_Z(V_1)=\kappa_0$,
it follows from \eqref{firstu} and \eqref{secondu} that
$\dim\Ss,\dim\hat\Ss\geq\kappa_0-\sigma+s_0$. Let 
$\delta\in\Ss\cap\hat\Ss$. Then 
$\delta=\pi_{\mathcal{V}(Z)}(\gamma_ZY_i)$ where $Y_i\in V_i$ and 
$\gamma_ZY_i\in\mathcal{R}$, $i=0,1$. By \eqref{disjoint} and the 
injectivity of $\pi_{\V(Z)}|_{\mathcal{R}}$ we have that
$\gamma_ZY_1=\gamma_ZY_0=0$. Thus $\delta=0$, and hence
$$
\dim\Ss\oplus\hat\Ss\geq 2(\kappa_0-\sigma+s_0).
$$
Since $r\leq s_0$, then that
$$
2(\kappa_0-\sigma+r)\leq 2(\kappa_0-\sigma+s_0)
\leq\dim\V^{p-\bar\tau,p-\bar\tau}(Z)
$$
concludes the proof of the claim.
\vspace{1ex}

Since $\mathcal{S}(\gamma)=W^{p,p}$ it holds that
\be\label{spanhatU0mlb}
\Sal(\hat\gamma)=\hat\Les^{\sigma}(X).
\ee
From \eqref{spanhatU0mlb} and $\sigma\geq 2$ we have  
$\hat\gamma\neq 0$. Since $\a$ is pluriharmonic then $\gamma$ 
symmetric and hence also is $\hat\gamma$. Thus \eqref{estpluri}, 
\eqref{partv} and \eqref{spanhatU0mlb} yield that
\be\label{kappaineqlmb}
4\sigma\leq\kappa_0(\kappa_0+2).
\ee

\noindent \emph{Case $\kappa_0=\sigma$}. This says 
that $\hat{\gamma}_Z(V)=\hat{\mathcal{L}}^{\sigma}(X)$ for
any $Z\in RE(\hat\gamma)$.
Given $Z\in RE(\hat\gamma)$ set 
$\gamma_1=\gamma_Z|_{N(X)}\colon N(X)\to\Les^{\sigma}(X)$
and $N_1=\ker\gamma_1$. Then $\dim N_1\geq\dim N(X)-\sigma$.
On one hand, if $\eta\in N(X)$ and $Y\in V^{2n}$ 
it follows from \eqref{decomp2mlb} that $\gamma_Y\eta=0$ 
if and only if $\lp\gamma_Y\eta,\hat{\gamma}_Z(V)\rp=0$. 
On the other hand, from  \eqref{decompmlb}, \eqref{decomp2mlb}
and the flatness of $\gamma$ we obtain
$$
\lp\gamma_Y\eta,\hat{\gamma}_Z(V)\rp
=\lp\gamma_Y\eta,\gamma_Z(V)\rp
=\lp\gamma_Y(V),\gamma_Z\eta\rp=0
$$
for any $\eta\in N_1$ and $Y\in V^{2n}$. Thus 
$N_1=\mathcal{N}(\gamma)$. Now \eqref{general} yields
\be\label{above2}
\nu(\gamma)=\dim N_1\geq\dim N(X)-\sigma
=2n-\kappa(\gamma)-\sigma
\ee
and gives \eqref{estim1mlb}.
\vspace{1ex}

We have seen that $\kappa_0$ and $2\leq\sigma\leq 10$ are 
both even. Since $\kappa_0\leq\sigma$ and the case of equality 
has already been considered then we assume that $\kappa_0<\sigma$.
Hence, in view of  \eqref{kappaineqlmb} it remains to consider the cases 
$(\kappa_0,\sigma)=(4,6),(6,8),(6,10)$ and $(8,10)$.
\vspace{1ex}

\noindent \emph{Cases $(6,8)$ and $(8,10)$}.
By \eqref{partv} the vector subspace 
$\hat\gamma_R(V)\cap\hat\gamma_S(V)$ is 
$\mathcal{T}|_{\Sal(\hat{\gamma})}$-invariant for any 
$R,S\in V^{2n}$ and thus of even dimension . Then by 
\eqref{spanhatU0mlb} there are $Z_1,Z_2\in RE(\hat\gamma)$ 
such that
\be\label{spanuchapelmlb}
\hat{\Les}^{\sigma}(X)=\hat{\gamma}_{Z_1}(V)+\hat{\gamma}_{Z_2}(V).
\ee

If $\eta\in N(X)$ and $Y\in V^{2n}$ it follows from 
\eqref{decomp2mlb} and \eqref{spanuchapelmlb} that
$\gamma_Y\eta=0$ if and only if
$\lp\gamma_Y\eta,\hat{\gamma}_{Z_j}(V)\rp=0$ for $j=1,2$.
Set $\gamma_1=\gamma_{Z_1}|_{N(X)}\colon N(X)\to\Les^{\sigma}(X)$,
$N_1=\ker\gamma_1$, 
$\gamma_2=\gamma_{Z_2}|_{N_1}\colon N_1\to\Les^{\sigma}(X)$
and $N_2=\ker\gamma_2$. Then $N_2=\mathcal{N}(\gamma)$ since 
from \eqref{decompmlb}, \eqref{decomp2mlb} and 
the flatness of $\gamma$ we have
$$
\lp\gamma_Y\eta,\hat{\gamma}_{Z_j}(V)\rp
=\lp\gamma_Y\eta,\gamma_{Z_j}(V)\rp\\
=\lp\gamma_Y(V),\gamma_{Z_j}\eta\rp=0,\;j=1,2,
$$
for any $\eta\in N_2$ and $Y\in V^{2n}$. From the Claim
above $\dim\gamma_{Z_j}(N(X))\leq 4$, $j=1,2$, and 
$$
\nu(\gamma)=\dim N_2\geq\dim N_1-4
\geq\dim N(X)-8\geq 2n-\kappa(\gamma)-\sigma
$$
as wished.
\vspace{1ex}

\noindent \emph{Case $(6,10)$}. 
If we have $Z_1,Z_2\in RE(\hat{\gamma})$ such that 
\eqref{spanuchapelmlb} holds then a similar argument as 
in the previous case gives \eqref{estim1mlb}. Otherwise, 
by \eqref{spanhatU0mlb} there are 
$Z_1,Z_2,Z_3\in RE(\hat{\gamma})$ such that
$$
\hat{\Les}^{10}(X)=\hat{\gamma}_{Z_1}(V)
+\hat{\gamma}_{Z_2}(V)+\hat{\gamma}_{Z_3}(V)
$$
and $\dim(\hat{\gamma}_{Z_1}(V)+\hat{\gamma}_{Z_2}(V))=8$.
Set $\gamma_1=\gamma_{Z_1}|_{N(X)}\colon N(X)\to\Les^{10}(X)$,
$N_1=\ker\gamma_1$,
\mbox{$\gamma_2=\gamma_{Z_2}|_{N_1}\colon N_1\to\Les^{10}(X)$},
$N_2=\ker\gamma_2$,
$\gamma_3=\gamma_{Z_3}|_{N_2}\colon N_2\to\Les^{10}(X)$
and $N_3=\ker\gamma_3$.
From \eqref{decompmlb} and \eqref{decomp2mlb} and the 
flatness of $\gamma$ we have
$$
\lp\gamma_{Z_3}\eta_2,\hat\gamma_{Z_j}(V)\rp
=\lp\gamma_{Z_3}\eta_2,\gamma_{Z_j}(V)\rp \\
=\lp\gamma_{Z_3}(V),\gamma_{Z_j}\eta_2)\rp=0,\,j=1,2, 
$$
for $\eta_2\in N_2$ and $Y\in V^{2n}$.  Hence
$\dim\gamma_{Z_3}(N_2)\leq 2$. Moreover, as in the previous 
case we obtain $N_3=\mathcal{N}(\gamma)$. From the Claim above 
we have $\dim\gamma_{Z_j}(N(X))\leq 4$, $j=1,2$, and
$$
\nu(\gamma)=\dim N_3\geq\dim N_2-2\geq\dim N_1-6
\geq\dim N(X)-10= 2n-\kappa(\gamma)-\sigma
$$
as wished.
\vspace{1ex}

\noindent \emph{Case $(4,6)$}. Given $Z_1\in RE(\hat\gamma)$ 
by \eqref{spanhatU0mlb} there is $Z_2\in RE(\hat\gamma)$ such 
that \eqref{spanuchapelmlb} holds. Suppose that there is 
$Z_1\in RE(\hat\gamma)$ such that $\dim\gamma_{Z_1}(N(X))\leq 4$. 
Since $\tau(\gamma)$ is even by \eqref{theclaim} this always 
holds if $\tau(\gamma)>6$. 
Set $\gamma_1=\gamma_{Z_1}|_{N(X)}\colon N(X)\to\Les^6(X)$ and 
$N_1=\ker\gamma_1$. 
From  \eqref{decompmlb} and \eqref{decomp2mlb} and the 
flatness of $\gamma$ we obtain
$$
\lp\gamma_{Z_2}\eta_1,\hat\gamma_{Z_1}(V)\rp
=\lp\gamma_{Z_2}\eta_1,\gamma_{Z_1}(V)\rp \\
=\lp\gamma_{Z_2}(V),\gamma_{Z_1}\eta_1)\rp=0
$$
for any $\eta_1\in N_1$.  Since $\kappa_0=4$ then
$\dim\gamma_{Z_2}(N_1)\leq 2$.

If $\eta\in N(X)$ and $Y\in V^{2n}$ it follows from 
\eqref{decomp2mlb} and \eqref{spanuchapelmlb} that
$\gamma_Y\eta=0$ if and only if
$\lp\gamma_Y\eta,\hat{\gamma}_{Z_j}(V)\rp=0$ for $j=1,2$. 
Set $\gamma_2=\gamma_{Z_2}|_{N_1}\colon N_1\to\Les^6(X)$ 
and $N_2=\ker\gamma_2$. As above we obtain that 
$N_2=\mathcal{N}(\gamma)$. Now \eqref{general} yields
$$
\nu(\gamma)=\dim N_2\geq\dim N_1-2
\geq\dim N(X)-6= 2n-\kappa(\gamma)-\sigma
$$
as wished.

By the above it remains to consider the case when 
$\tau(\gamma)=6$ and 
\be\label{aboveas}
\gamma_Z(N(X))=\Les^6(X)
\;\;\mbox{for any}\;\; Z\in RE(\hat{\gamma}).
\ee
If $Y\in RE(\gamma)$ then $\sigma_\gamma(Y)\leq\tau(\gamma)=6$ 
by \eqref{sigmatau}. Suppose that there is $Y\in RE(\gamma)$ 
such that $\sigma_\gamma(Y)\leq 4$. From \eqref{kappaineqlmb} 
we are in case $\kappa_0=\sigma$ for $Y$ and thus 
$\nu(\gamma)\geq 2n-k(\gamma)-\sigma_\gamma(Y)$.
Since $\sigma_\gamma(Y)<6=\sigma$ then \eqref{estim1mlb} 
also holds for $X$.  

In view of the above, we assume further that 
$\sigma_\gamma(Y)=6$ for any $Y\in RE(\gamma)$.
Now let $Z_1\in RE(\gamma)\cap RE(\hat\gamma)$ and then let 
$\tilde{\gamma}\colon V^{2n}\times V^{2n}\to\hat{\mathcal{L}}^6(Z_1)$ 
stand for taking the $\hat{\mathcal{L}}^6(Z_1)$-component of 
$\gamma$. Suppose that there is $Z_2\in RE(\tilde\gamma)$ such 
that $\tilde{\gamma}_{Z_2}(V)=\hat{\mathcal{L}}^6(Z_1)$. 
Under this assumption for $Z_1$ we  are in the situation analyzed 
in the Case $\kappa_0=\sigma$ and thus \eqref{estim1mlb} holds for 
$Z_1$. Since $\sigma_\gamma(Z_1)=\sigma$ it also holds for $X$. 

In view of \eqref{kappaineqlmb} we now also assume that 
$\dim\tilde{\gamma}_{Z_2}(V)=4$ for any $Z_2\in RE(\tilde\gamma)$.
If $\dim\gamma_{Z_2}(N(Z_1))\leq 4$ for some $Z_2\in RE(\tilde\gamma)$ 
then the initial part of the proof of this case gives that 
\eqref{estim1mlb} holds for $Z_1$ and then also for $X$ since 
$\sigma_\gamma(Z_1)=\sigma$. Hence we assume that 
$\gamma_{Z_2}(N(Z_1))=\mathcal{L}(Z_1)$ for any 
$Z_2\in RE(\tilde\gamma)$.  

The remaining case to consider is when there are 
$Z_1,Z_2\in RE(\gamma)\cap RE(\hat\gamma)$ and $Z_2\in RE(\tilde\gamma)$
for which \eqref{spanuchapelmlb} holds, $\sigma_\gamma(Z_j)=6$, $j=1,2$, 
$\gamma_{Z_2}(N(Z_1))=\Les(Z_1)$ and  $\dim\tilde\gamma_{Z_2}(V)=4$. 
To conclude the proof we show that this situation is not possible. 
Hence suppose otherwise. In particular, we have 
$\Les(Z_1)\subset\gamma_{Z_2}(V)$. 
From \eqref{aboveas} we obtain that $\Les(X)\subset\gamma_{Z_j}(V)$, $j=1,2$. 
Thus given $\eta_0\in\Les(X)$ there are $Y_1,Y_2\in V^{2n}$ such that 
$\eta_0=\gamma_{Z_1}Y_1=\gamma_{Z_2}Y_2$.  
Let $\xi_0\in\Les(X)$ and $\xi_j\in\Les(Z_j)$, $j=1,2$. Then 
$$
\lp\xi_0+\xi_1+\xi_2,\eta_0\rp
=\lp\xi_1,\gamma_{Z_1}Y_1\rp+\lp\xi_2,\gamma_{Z_2}Y_2\rp=0.
$$
If $\eta_1\in\mathcal{L}(Z_1)$ then $\lp\xi_0,\eta_1\rp=0$ 
since $\mathcal{L}(X)\subset\gamma_{Z_1}(V)$ and
$\mathcal{L}(Z_1)\subset\U(Z_1)$. Let $Y_3\in V^{2n}$
be such that $\eta_1=\gamma_{Z_2}Y_3$.  Since
$\mathcal{L}(Z_2)\subset\U(Z_2)$ then
$$
\lp\xi_0+\xi_1+\xi_2,\eta_1\rp=\lp\xi_2,\gamma_{Z_2}Y_3\rp=0.
$$
If $\eta_2\in\mathcal{L}(Z_2)$ then 
$\lp\xi_j,\eta_2\rp=0$, $j=0,1$, since 
$\mathcal{L}(X)\subset\gamma_{Z_2}(V)$, 
$\mathcal{L}(Z_1)\subset\gamma_{Z_2}(V)$ and
$\mathcal{L}(Z_2)\subset\U(Z_2)$.  
Thus $\lp\xi_0+\xi_1+\xi_2,\eta_2\rp=0$.  
Hence $\Les(X)+\Les(Z_1)+\Les(Z_2)$ is an isotropic
vector subspace.

We argue that $\dim\Les(X)\cap\Les(Z_j)=\dim\Les(Z_1)\cap\Les(Z_2)=2$.
On one hand, we have that $\Les(X)\cap\Les(Z_j)\neq 0$ since otherwise
the vector subspace $\Les(X)\oplus\Les(Z_j)$ would be isotropic of 
dimension $12$ which is not possible. On the other hand, we have
$$
\lp\xi,\hat\gamma_{Z_j}(V)\rp=\lp\xi,\gamma_{Z_j}(V)\rp=0
$$
for any $\xi\in\Les(X)\cap\Les(Z_j)$. Since $\kappa_0=4$ it
follows from part $(ii)$ of Proposition \ref{even} that
$\dim\Les(X)\cap\Les(Z_j)=2$. Having that
$\mathcal{L}(Z_1)\oplus\mathcal{L}(Z_2)\subset\gamma_{Z_2}(V)$ 
is isotropic yields $\mathcal{L}(Z_1)\cap\mathcal{L}(Z_2)\neq 0$. 
If $\xi\in\mathcal{L}(Z_1)\cap\mathcal{L}(Z_2)$ then 
$$
\lp\xi,\tilde\gamma_{Z_2}(V)\rp=\lp\xi,\gamma_{Z_2}(V)\rp=0
$$
where the second equality follows from
$\mathcal{L}(Z_2)\subset\U(Z_2)$. Since 
$\dim\tilde\gamma_{Z_2}(V)=4$ then 
$\mathcal{L}(Z_1)\cap\mathcal{L}(Z_2)=2$.  
We have shown that $\Les(X)+\Les(Z_1)+\Les(Z_2)$ has dimension
$12$, but this is a contradiction.\qed 

\begin{remark}\po {\em The estimate 
$\nu(\gamma)\geq 2n-2p$ is Proposition $10$ in \cite{CCD}.
A counterexample constructed in \cite{CCD} 
shows that this estimate is false already for $p=12$.
}\end{remark}

Henceforward $U^p=U_1^s\oplus U_2^{p-s}$ is an orthogonal 
decomposition where 
$$
U_1^s=\Sal(\pi_1\circ\beta).
$$

\begin{lemma}\po\label{sbeta} If \eqref{productflat2} holds then
\be\label{betaimage1} 
\mathcal{S}(\beta)=U_1^s\oplus U_1^s
\ee
and $\mathcal{N}(\beta)=\mathcal{N}(\gamma_{U_1})$ 
where $\gamma_{U_1}=\pi_{U_1\times U_1}\circ\gamma$.
\end{lemma}

\proof We have that
$$
\beta(X,Y)=(\xi,\eta)\iff
\beta(Y,X)=(\xi,-\eta)\iff\beta(X,JY)=(\eta,-\xi).
$$ 
Thus  if $(\xi,\eta)=\sum_k\beta(X_k,Y_k)$ then
$$
\sum_k\beta(Y_k,X_k)=(\xi,-\eta),\;\;\sum_k\beta(X_k,JY_k)
=(\eta,-\xi),\;\;\sum_k\beta(JY_k,X_k)
=(\eta,\xi).
$$
Hence if $(\xi,\eta)\in\Sal(\beta)$ then 
$(\xi,0),(0,\xi),(\eta,0)\in\Sal(\beta)$ and thus 
$\mathcal{S}(\beta)\subset U_1\oplus U_1$. On the other 
hand, if $(\xi,\eta)\in U_1\oplus U_1$ there are 
$\bar{\xi},\bar{\eta}\in U^p$ so that 
$(\xi,\bar{\xi}),(\eta,\bar{\eta})\in\Sal(\beta)$
and thus $(\xi,\eta)\in\Sal(\beta)$ which proves 
\eqref{betaimage1}. 
\vspace{1ex}

From \eqref{productflat2}, \eqref{betaimage1} and 
$(U_2\oplus U_2)^\perp=U_1\oplus U_1$ we   
obtain $\mathcal{S}(\gamma|_{V\times\mathcal{N}(\beta)})
\subset U_2\oplus U_2$. Then
$\lp\gamma(X,Y),(\xi,\bar\xi)\rp=0$ if $X\in V^{2n}$,
$Y\in\mathcal{N}(\beta)$ and $\xi,\bar\xi\in U_1^s$. 
Thus $\mathcal{N}(\beta)\subset\mathcal{N}(\gamma_{U_1})$.  
On the other hand, if $S\in\mathcal{N}(\gamma_{U_1})$ 
then $\beta(X,S)=\gamma_{U_2}(X,S)+\gamma_{U_2}(JX,JS)$=0
by \eqref{betaimage1} for any $X\in V^{2n}$. \qed

\begin{lemma}\po\label{diagonalization}
Let the bilinear form $\beta\colon V^{2n}\times V^{2n}\to W^{p,p}$, 
$p\leq n$, be flat. Then 
\be\label{maincostum}
\nu(\beta)=2n-\kappa(\beta). 
\ee
Moreover, if $\kappa(\beta)=2p$ there is a 
basis $\{X_i,JX_i\}_{1\leq i\leq n}$ of $V^{2n}$ such that:
\begin{itemize}
\item[(i)] $\mathcal{N}(\beta)=\spa\{X_j,JX_j,\; p+1\leq j\leq n\}$.
\item[(ii)] $\beta(Y_i,Y_j)=0\;\;\mbox{if}\;\;i\neq j$ 
and $Y_k\in\spa\{X_k,JX_k\}\;\;\mbox{for}\;\;k=i,j$. 
\item[(iii)] $\{\beta(X_j,X_j),\beta(X_j,JX_j)\}_{1\leq j\leq p}$ 
is an orthonormal basis of $W^{p,p}$.
\end{itemize}
\end{lemma}

\proof Proposition $7$ in \cite{CCD} gives \eqref{maincostum}.
The remaining of the statement is Proposition~$2.6$ in 
\cite{CDa} as well as Lemma $7$ in \cite{FHZ}.\vspace{2ex}\qed

Let $\theta\colon V^{2n}\times V^{2n}\to W^{p,p}$ be 
the pluriharmonic symmetric bilinear form defined as
\be\label{theta}
\theta(X,Y)=\gamma(X,Y)-\gamma(JX,JY).
\ee 
It follows from \eqref{T1} that
\be\label{T3}
\mathcal{T}\theta(X,Y)=\theta(X,JY)\;\;\mbox{for any}\;\;X,Y\in V^{2n}.
\ee

If the condition \eqref{productflat2} holds then also
\be\label{productflat}
\lp\beta(X,Y),\theta(Z,T)\rp=\lp\beta(X,T),\theta(Z,Y)\rp\;\;
\mbox{for any}\;\;X,Y,Z,T\in V^{2n}.
\ee
In fact, we have using \eqref{T1} and \eqref{T2} that
\begin{align*}
\lp\beta(X,Y),\theta(Z,T)\rp
&=\lp\beta(X,Y),\gamma(Z,T)\rp-\lp\beta(X,Y),\gamma(JZ,JT)\rp\\
&=\lp\beta(X,T),\gamma(Z,Y)\rp-\lp\beta(X,JT),\gamma(JZ,Y)\rp\\
&=\lp\beta(X,T),\gamma(Z,Y)\rp-\lp\mathcal{T}\beta(X,T),\gamma(JZ,Y)\rp\\
&=\lp\beta(X,T),\gamma(Z,Y)\rp-\lp\beta(X,T),\gamma(JZ,JY)\rp\\
&=\lp\beta(X,T),\theta(Z,Y)\rp.
\end{align*}

If $\gamma$ is flat then also is $\theta$. 
In effect, using \eqref{T1} we have
\begin{align*}
\lp\theta(X,Y)&,\theta(Z,T)\rp
=\lp\gamma(X,Y),\gamma(Z,T)\rp-\lp\gamma(X,Y),\gamma(JZ,JT)\rp\\
&\quad-\lp\gamma(JX,JY),\gamma(Z,T)\rp+\lp\gamma(JX,JY),\gamma(JZ,JT)\rp\\
&=\lp\gamma(X,T),\gamma(Z,Y)\rp-\lp\gamma(X,JT),\gamma(JZ,Y)\rp\\
&\quad-\lp\gamma(JX,T),\gamma(Z,JY)\rp+\lp\gamma(JX,JT),\gamma(JZ,JY)\rp\\
&=\lp\gamma(X,T),\gamma(Z,Y)\rp-\lp\gamma(X,T),\gamma(JZ,JY)\rp\\
&\quad-\lp\gamma(JX,JT),\gamma(Z,Y)\rp+\lp\gamma(JX,JT),\gamma(JZ,JY)\rp\\
&=\lp\theta(X,T),\theta(Z,Y)\rp.
\end{align*}

Since $2\gamma=\beta+\theta$ then 
$\mathcal{N}(\beta)\cap\mathcal{N}(\theta)\subset\mathcal{N}(\gamma)$
whereas the opposite inclusion follows from \eqref{beta}, 
\eqref{theta} and that $\mathcal{N}(\gamma)$ is $J$-invariant. 
Therefore
\be\label{nucleogamma}
\mathcal{N}(\gamma)=\mathcal{N}(\beta)\cap\mathcal{N}(\theta)
\ee
and, in particular, we have that $\nu(\beta)\geq\nu(\gamma)$. 

\begin{lemma}\po\label{bgbasicfacts} Let 
$\gamma,\beta\colon V^{2n}\times V^{2n}\to W^{p,p}$, $p\leq n$, 
be flat and satisfy the condition \eqref{productflat2}. 
If $\nu(\beta)=2n-2s$ then the bilinear forms 
$\theta_j=\pi_{U_j\times U_j}\circ\theta$, 
$j=1,2$, are flat.
\end{lemma}

\proof By \eqref{betaimage1} and 
Lemma \ref{diagonalization} there is a basis 
$\{X_j,JX_j\}_{1\leq j\leq n}$ of $V^{2n}$ satisfying that 
$\mathcal{N}(\beta)=\spa\{X_j,JX_j,\;s+1\leq j\leq n\}$, that
\be\label{diagonal}
\beta(X_i,X_j)=0=\beta(X_i,JX_j)\;\;\mbox{if}\;\;i\neq j 
\ee 
and that
$\{\beta(X_j,X_j),\beta(X_j,JX_j)\}_{1\leq j\leq s}$
is an orthonormal basis of $U_1^s\oplus U_1^s$.

From \eqref{betaimage1} we have that \eqref{productflat} is 
equivalent to
\be\label{flatbeta}
\lp\beta(X,Y),\theta_1(Z,T)\rp=\lp\beta(X,T),\theta_1(Z,Y)\rp
\;\;\mbox{for any}\;\; X,Y,Z,T\in V^{2n}.
\ee
In particular, we obtain using \eqref{diagonal} that
$$
\theta_1(X_i,X_i),\theta_1(X_i,JX_i)\in
\spa\{\beta(X_i,X_i),\beta(X_i,JX_i)\}.
$$ 
Moreover, since $\theta_1$ is symmetric it follows from 
\eqref{diagonal} and \eqref{flatbeta} for $k\neq\ell$ that
$$
\lp\beta(X_j,X_j),\theta_1(X_k,X_\ell)\rp
=\lp\beta(X_j,X_\ell),\theta_1(X_k,X_j)\rp
=\lp\beta(X_j,X_k),\theta_1(X_\ell,X_j)\rp=0,
$$
$$
\!\!\lp\beta(X_j,JX_j),\theta_1(X_k,X_\ell)\rp
=\lp\beta(X_j,X_\ell),\theta_1(X_k,JX_j)\rp
=\lp\beta(X_j,X_k),\theta_1(X_\ell,JX_j)\rp=0
$$
and thus
$$
\theta_1(X_k,X_\ell)=0=\theta_1(X_k,JX_\ell)
\;\;\mbox{if}\;\; k\neq\ell.
$$ 
We have that
$$
\mathcal{T}\gamma_{U_1}(X,Y)=\mathcal{T}(\a_{U_1}(X,Y),\a_{U_1}(X,JY))=
(\a_{U_1}(X,JY),-\a_{U_1}(X,Y))=\gamma_{U_1}(X,JY)
$$
for any $X,Y\in V^{2n}$. Then 
$$
\mathcal{T}\theta_1(X,Y)
=\mathcal{T}(\gamma_{U_1}(X,Y)-\gamma_{U_1}(JX,JY))
=\gamma_{U_1}(X,JY)+\gamma_{U_1}(JX,Y)
=\theta_1(X,JY)
$$
and hence
$$
\lp\theta_1(X_i,X_i),\theta_1(JX_i,JX_i)\rp
=\lp\theta_1(X_i,JX_i),\theta_1(X_i,JX_i)\rp.
$$
We have shown that $\theta_1$ and $\theta$ are flat. Since 
$\theta=\theta_1\oplus\theta_2$ then also $\theta_2$ is flat.\qed

\begin{lemma}\po
Let $\gamma,\beta\colon V^{2n}\times V^{2n}\to W^{p,p}$,
$p\leq n$, be flat and satisfy the condition~\eqref{productflat2}.
If $\Sal(\gamma)=W^{p,p}$ and $p-s\leq 9$ 
the flat bilinear form 
$\varphi=\theta|_{V\times\mathcal{N}(\beta)}
\colon V^{2n}\times\mathcal{N}(\beta)\to W^{p,p}$ 
satisfies
\be\label{bnewmainlemma2}
0\leq\nu(\beta)-\nu(\gamma)\leq\kappa(\varphi)
+\sigma(\varphi)\leq 2p-2s.
\ee
\end{lemma}

\proof It follows from \eqref{betaimage1} and \eqref{productflat} 
that $\Sal(\varphi)\subset U_2^{p-s}\oplus U_2^{p-s}\subset W^{p,p}$.
Thus $\varphi$ is seen  in the sequel as a map  
\be\label{deffi}
\varphi\colon V^{2n}\times\mathcal{N}(\beta)
\to U_2^{p-s}\oplus U_2^{p-s}.
\ee
To obtain the proof it suffices to show for any $X\in RE(\varphi)$ 
that we have 
\be\label{l15inter}
0\leq\nu(\beta)-\nu(\gamma)\leq\kappa(\varphi)+\sigma_\varphi(X)
\leq 2p-2s.
\ee 
Fix $X\in RE(\varphi)$ and set $\sigma=\sigma_\varphi(X)$ for
simplicity. Proposition \ref{decomp} gives a decomposition 
\be\label{decomp2}
U_2^{p-s}\oplus U_2^{p-s}
=\U^\tau(X)\oplus\hat{\U}^\tau(X)\oplus\V^{p-s-\tau,p-s-\tau}
\ee
where  $\U^\tau(X)=\varphi_X(\mathcal{N}(\beta))\cap
\varphi_X(\mathcal{N}(\beta))^\perp$ with
$\tau=\tau_\varphi(X)$ for simplicity of notation and 
$\varphi_X(\mathcal{N}(\beta))\subset\U^\tau(X)
\oplus\V^{p-s-\tau,p-s-\tau}$. 
Thus $\kappa(\varphi)\leq 2p-2s-\tau$. From \eqref{secondst}
we obtain that $\sigma(\varphi)\leq\sigma\leq\tau$  
and the last inequality in \eqref{l15inter} follows.
\vspace{1ex}

We have that $0\leq\sigma\leq\tau\leq p-s\leq 9$. From
\eqref{T3} and the  $J$-invariance of $\mathcal{N}(\beta)$ 
we obtain that $\mathcal{T}\varphi(Y,Z)=\varphi(Y,JZ)$ for 
any $Y\in V^{2n}$ and $Z\in\mathcal{N}(\beta)$. Thus 
part $(iii)$ of Proposition \ref{even} gives that 
$N(X)=\ker\varphi_X$ is $J$-invariant and hence part 
$(ii)$ that $\sigma$ is even. Therefore, we have that
$0\leq\sigma\leq 8$.
\vspace{1ex}

\noindent\emph{Case $\sigma=0$}.
Since $\theta(Y,N(X))=0$ for any $Y\in V^{2n}$ then
$N(X)\subset\mathcal{N}(\beta)\cap\mathcal{N}(\theta)
=\mathcal{N}(\gamma)$ from \eqref{nucleogamma}.
On the other hand, it is a general fact that
\be\label{estimateN(X)}
\dim N(X)=\nu(\beta)-\kappa(\varphi)
\ee
and since $\sigma=0$ then \eqref{l15inter} follows 
from \eqref{estimateN(X)}.
 
Henceforward, we assume that $\sigma\geq 2$. As in 
\eqref{decomp2mlb} we have a decomposition
\be\label{decomp0}
\U^\tau(X)\oplus\hat\U^\tau(X)
=\Les^{\sigma}(X)\oplus\hat\Les^{\sigma}(X)
\oplus\V_0^{\tau-\sigma,\tau-\sigma}
\ee
We claim that the symmetric bilinear form defined by 
$\hat\theta=\pi_{\hat\Les^{\sigma}(X)}\circ\theta$ satisfies
\be\label{spanhatU0}
\Sal(\hat\theta)=\hat\Les^{\sigma}(X).
\ee
In fact, if otherwise there is $0\neq\eta
=\sum_{j=1}^r\varphi(Y_j,S_j)\in\mathcal{L}^{\sigma}(X)$
with $Y_1,\ldots,Y_r\in V^{2n}$ and $S_1,\ldots,S_r\in N(X)$ 
such that
\be\label{prevs}
0=\lp\eta,\hat\theta(Z,T)\rp
=\sum_{j=1}^r\lp\theta(Y_j,S_j),\theta(Z,T)\rp
\ee
for any $Z,T\in V^{2n}$.  Since $2\gamma=\beta+\theta$ 
we obtain from \eqref{productflat} and \eqref{prevs} that
$$
2\lp\eta,\gamma(Z,T)\rp
=\sum_{j=1}^r\lp\theta(Y_j,S_j),\beta(Z,T)+\theta(Z,T)\rp
=\sum_{j=1}^r\lp\theta(Y_j,S_j),\theta(Z,T)\rp=0
$$
for any $Z,T\in V^{2n}$. Since $0\neq\eta=2\sum_{j=1}^r\gamma(Y_j,S_j)$ 
this is a contradiction by the assumption for $\Sal(\gamma)$  
that proves the claim.

Part $(i)$ of Proposition \ref{even} gives that
$\mathcal{T}|_{\mathcal{L}^{\sigma}(X)}
\in\text{Aut}(\mathcal{L}^{\sigma}(X))$ 
is a complex structure. Then from \eqref{decomp2} and \eqref{decomp0} 
we obtain for any $(\xi,\bar\xi)\in\mathcal{L}^{\sigma}(X)$ that
$$
\lp\mathcal{T}\hat\theta(Z,Y),(\xi,\bar\xi)\rp
=\lp\theta(Z,Y),\mathcal{T}(\xi,\bar\xi)\rp
=\lp\theta(Z,JY),(\xi,\bar\xi)\rp
=\lp\hat\theta(Z,JY),(\xi,\bar\xi)\rp
$$
which says that
\be\label{welldef}
\mathcal{T}\hat{\theta}{}_ZY
=\hat{\theta}{}_ZJY
\ee
for any $Z,Y\in V^{2n}$. Hence if $Z\in RE(\hat\theta)$ then 
$\kappa_0=\kappa(\hat\theta)$ is even.  Being $\theta$ symmetric
then also is $\hat\theta$ and it follows from 
\eqref{estpluri} and \eqref{spanhatU0} that
\be\label{kappaineq}
4\sigma\leq\kappa_0(\kappa_0+2).
\ee
Therefore, if  $\sigma=2,4$ then $\sigma=\kappa_0$ 
and if $\sigma=6,8$ then either $\sigma=\kappa_0$ 
or $\sigma=\kappa_0+2$.
\vspace{1ex}

\noindent\emph{Case $\sigma=\kappa_0$}. Given $Z\in RE(\hat\theta)$
we have $\hat\theta_Z(V)=\hat\Les^{\sigma}(X)$. Since
$\varphi(Y,\eta)\in\Les^{\sigma}(X)$ if $Y\in V^{2n}$ and 
$\eta\in N(X)$, then
\be\label{condition}
\varphi(Y,\eta)=0\;\;\mbox{if and only if}
\;\;\lp\varphi(Y,\eta),\hat\theta_ZT\rp=0
\ee
for any $T\in V^{2n}$.
Set $\theta_1=\theta_Z|_{N(X)}\colon N(X)\to\Les^{\sigma}(X)$ and
$N_1=\ker\theta_1$. From \eqref{decomp2}, \eqref{decomp0} and the 
flatness of $\theta$ we obtain
$$
\lp\varphi(Y,\delta),\hat\theta_ZT\rp
=\lp\theta(Y,\delta),\theta(Z,T)\rp
=\lp\theta(Y,T),\theta(Z,\delta)\rp=0
$$
for any $\delta\in N_1$ and $Y,T\in V^{2n}$.  Hence from 
\eqref{condition} we have 
$N_1\subset\mathcal{N}(\beta)\cap\mathcal{N}(\theta)
=\mathcal{N}(\gamma)$. Then \eqref{nucleogamma} and
\eqref{estimateN(X)} give 
$$
\nu(\gamma)\geq \dim N_1\geq \dim N(X)-\sigma=
\nu(\beta)-\kappa(\varphi)-\sigma,
$$
and \eqref{l15inter} follows. 

\noindent\emph{Case $\sigma=\kappa_0+2$}.
Suppose that there is $Z\in RE(\varphi)$ such that
$\Les(Z)=\Sal(\varphi|_{V\times N(Z)})$ satisfies
$\sigma_\varphi(Z)\leq 4$. Then by \eqref{kappaineq} 
for such $Z\in RE(\varphi)$ we are in Case $\sigma=\kappa_0$. 
Hence
$$
\nu(\beta)-\nu(\gamma)\leq\kappa(\varphi)
+\sigma_\varphi(Z)\leq \kappa(\varphi)+4<\kappa(\varphi)+\sigma,
$$
and since $\sigma\geq 6$ then \eqref{l15inter} holds. 
Thus henceforward we assume that 
$\sigma_\varphi(Z)\geq 6$ for any $Z\in RE(\varphi)$.

If $Z_1,Z_2\in RE(\hat\theta)$ then \eqref{welldef} gives 
that $\hat\theta_{Z_1}(V)\cap\hat\theta_{Z_2}(V)$ is 
$\mathcal{T}$-invariant and therefore of even dimension.
Given $Z_1\in RE(\hat\theta)$ then by \eqref{spanhatU0} there is
$Z_2\in RE(\hat\theta)$ such that we have
$\hat\Les^{\sigma}(X)=\hat\theta_{Z_1}(V)+\hat\theta_{Z_2}(V)$.

Since $\varphi(Y,\eta)\in\Les^{\sigma}(X)$
if $\eta\in N(X)$ and $Y\in V^{2n}$, then
\be\label{condition2}
\varphi(Y,\eta)=0\;\;\mbox{if and only if}\;\;
\lp\varphi(Y,\eta),\hat\theta_{Z_j}T\rp=0,\; j=1,2,
\ee
for any $T\in V^{2n}$. If 
$\theta_1=\theta_{Z_1}|_{N(X)}\colon N(X)\to\Les^{\sigma}(X)$ 
and $N_1=\ker\theta_1$, then
\be\label{estimateN_1}
\dim N(X)=\dim N_1+\dim\theta_1(N(X)).
\ee
If $\theta_2=\theta_{Z_2}|_{N_1}\colon N_1\to\Les^{\sigma}(X)$
and $N_2=\ker\theta_2$ we have
from \eqref{decomp2} and \eqref{decomp0} that
$$
\lp\theta_2\delta_1,\hat\theta_{Z_1}Y\rp
=\lp\theta(Z_2,\delta_1),\theta_{Z_1}Y\rp
=\lp\theta(Z_2,Y),\theta(Z_1,\delta_1)\rp=0
$$
for any $\delta_1\in N_1$ and $Y\in V^{2n}$.  Thus
$\dim\theta_2(N_1)\leq\sigma-\kappa_0=2$ and hence
\be\label{estimateN_2}
\dim N_1\leq\dim N_2+2.
\ee
It follows from \eqref{decomp2} and \eqref{decomp0} that
$$
\lp\varphi(Y,\delta_2),\hat\theta_{Z_j}T\rp
=\lp\theta(Y,\delta_2),\theta(Z_j,T)\rp
=\lp\theta(Y,T),\theta(Z_j,\delta_2)\rp=0
$$
for any $\delta_2\in N_2$, $Y,T\in V^{2n}$ and $j=1,2$.
Then from \eqref{condition2} we have
$N_2\subset\mathcal{N}(\beta)\cap\mathcal{N}(\theta)$.
Hence using \eqref{nucleogamma}, \eqref{estimateN(X)}, 
\eqref{estimateN_1} and \eqref{estimateN_2} we obtain
\begin{align*}
\nu(\gamma)\geq&\dim N_2\geq\dim N_1-2
=\dim N(X)-\dim\theta_1(N(X))-2\\
&=\nu(\beta)-\kappa(\varphi)-\dim\theta_1(N(X))-2.
\end{align*}

Since $\sigma=\kappa_0+2$ then in order from the above to 
have \eqref{l15inter} it is necessary to show that there 
is $Z\in RE(\hat\theta)$ such that $\dim\theta_{Z}(N(X))\leq\kappa_0$. 
Arguing by contradiction, assume 
that $\dim\theta_{Z}(N(X))>\kappa_0$ for any $Z\in RE(\hat\theta)$.
Since $N(X)$ is $J$-invariant then $\theta_{Z}(N(X))$
has even dimension and thus the assumption means that
$\theta_{Z}(N(X))=\Les^{\sigma}(X)$ for any $Z\in RE(\hat\theta)$.
Let us take $Z\in RE(\varphi)\cap RE(\theta)\cap RE(\hat\theta)$.
Then $N(Z)=\ker\varphi_{Z}$ satisfies $N(Z)\subset\ker\theta_Z$. 
From \eqref{secondst} and \eqref{deffi} we have 
\be\label{cal}
\Les^{\sigma_\varphi(Z)}(Z)
\subset{\cal R}=\theta_Z(V)\cap\theta_{Z}(V)^\perp
\cap(U_2^{p-s}\oplus U_2^{p-s}) 
\ee
and from \eqref{decomp2} and \eqref{decomp0} that there is 
a decomposition
$$
U_2^{p-s}\oplus U_2^{p-s}=\Les^{\sigma}(X)\oplus\hat\Les^{\sigma}(X)
\oplus\V_0^{\tau-\sigma,\tau-\sigma}
\oplus\V^{p-s-\tau,p-s-\tau}.
$$
Thus if $\delta\in{\cal R}$ then $\delta=\delta_1+\delta_2+\delta_3$
where $\delta_1\in\Les^{\sigma}(X)$, $\delta_2\in\hat\Les^{\sigma}(X)$
and $\delta_3\in\V_0\oplus\V$.
Since $\delta\in\theta_{Z}(V)^\perp$ and  
$\theta_{Z}(N(X))=\Les^{\sigma}(X)$ then
$$
0=\lp\delta,\theta_Z(N(X))\rp=\lp\delta_2,\Les^{\sigma}(X)\rp
$$
and hence $\delta_2=0$. Thus if $\delta^1,\delta^2\in{\cal R}$ 
we have
$$
0=\lp\delta^1,\delta^2\rp
=\lp\delta^1_1+\delta^1_3,\delta^2_1+\delta^2_3\rp
=\lp\delta^1_3,\delta^2_3\rp.
$$
Hence if $\pi\colon{\cal R}\to\V_0\oplus\V$ is
defined by $\pi(\delta)=\delta_3$, then the vector subspace
$\pi({\cal R})\subset\V_0\oplus\V$ is isotropic and thus 
$\dim\pi({\cal R})\leq p-s-\sigma$. Since $p-s\leq 9$ 
by assumption and $\sigma\geq 6$ then $\dim\pi({\cal R})\leq 3$. 
From \eqref{cal} we have $\dim{\cal R}\geq\sigma_\varphi(Z)$, 
and hence $\dim\ker\pi\geq\sigma_\varphi(Z)-3\geq 3$. 
On the other hand, we have that 
$\ker\pi\subset\Les^{\sigma}(X)\cap{\cal R}$.  
Since ${\cal R}\subset\theta_Z(V)^\perp$, we obtain from 
\eqref{decomp2} and \eqref{decomp0} that
$$
\lp\zeta,\hat\theta_{Z}Y\rp=\lp\zeta,\theta_{Z}Y\rp=0
$$
for any $\zeta\in\ker\pi$ and $Y\in V^{2n}$. Therefore,  
that $\ker\pi\subset\Les^{\sigma}(X)$,
$\hat\theta_{Z}(V)\subset\hat\Les^{\sigma}(X)$ and that
$\dim\hat\theta_{Z}(V)=\kappa_0=\sigma-2$ yield
$\dim\ker\pi\leq 2$, and we reached a contradiction.\qed

\begin{lemma}\po\label{alglemma} Let
$\gamma,\beta\colon V^{2n}\times V^{2n}\to W^{p,p}$, $p\leq n$, 
be flat and satisfy the condition \eqref{productflat2}. If 
$p\leq 11$ and $\Sal(\gamma)=W^{p,p}$ then 
$\nu(\gamma)\geq 2n-2p$.
\end{lemma}

\proof First assume that $s\geq 2$ in which case 
\eqref{bnewmainlemma2} holds. Let  
$\varphi\colon V\times\mathcal{N}(\beta)\to U_2^{p-s}\oplus U_2^{p-s}$
be given by \eqref{deffi} and let $X\in RE(\varphi)$ 
satisfy $\sigma_\varphi(X)=\sigma(\varphi)$.  
By \eqref{decomp2} we have that
$$
U_2^{p-s}\oplus U_2^{p-s}=\U^\tau(X)\oplus\hat\U^\tau(X)
\oplus\V^{p-s-\tau,p-s-\tau}
$$
with $\varphi_X(\mathcal{N}(\beta))\subset\U^\tau(X)
\oplus\V^{p-s-\tau,p-s-\tau}$. Then
$\kappa(\varphi)\leq 2p-2s-\tau$.
On the other hand, it follows from \eqref{betaimage1} and 
\eqref{maincostum} that $\nu(\beta)\geq 2n-2s$ whereas by 
\eqref{bnewmainlemma2} we have that 
$\nu(\beta)-\nu(\gamma)\leq\kappa(\varphi)+\sigma(\varphi)$. 
Hence
$$
2n-2s-\nu(\gamma)\leq\nu(\beta)-\nu(\gamma)
\leq\kappa(\varphi)+\sigma(\varphi)
\leq 2p-2s-\tau+\sigma(\varphi).
$$
Since $\tau\geq\sigma(\varphi)$ by \eqref{secondst} then 
$\nu(\gamma)\geq 2n-2p$ as we wished.
\vspace{1ex}

If $s=0$ then $\beta=0$, that is, $\a$ is pluriharmonic and 
then Lemma \ref{mainlemmabis} gives the result. Thus it remains 
to consider the case $s=1$ and thus $\beta\neq 0$. Part $(ii)$ 
of Proposition \ref{even} yields that the vector space 
$\beta_X(V)$ is even dimensional. 
Thus $\kappa(\beta)=2$ and then \eqref{maincostum} gives that 
$\nu(\beta)=2n-2$. Then by Lemma \ref{bgbasicfacts} we have that
$\theta=\theta_1+\theta_2$ where the bilinear form
$\theta_2\colon V^{2n}\times V^{2n}\to U_2^{p-1}\oplus U_2^{p-1}$ 
is flat.
 
We claim that $\theta_2$ satisfies the assumptions of 
Lemma \ref{mainlemmabis}. From \eqref{betaimage1} and 
$2\gamma=\beta+\theta$ we have 
$$
\theta_2(X,Y)=2\gamma_2(X,Y)=2(\a_{U_2}(X,Y),\a_{U_2}(X,JY)).
$$
Since $\Sal(\gamma)=W^{p,p}$ by assumption then 
$\Sal(\theta_2)=U_2^{p-1}\oplus U_2^{p-1}$.
If follows from \eqref{betaimage1} that the symmetric bilinear 
form $\a_{U_2}$ is pluriharmonic and hence also is $\theta_2$. 

By Lemma \ref{mainlemmabis} we have $\nu(\theta_2)\geq 2n-2p+2$.
Then \eqref{theta} yields 
$\mathcal{N}(\gamma_{U_1})\subset\mathcal{N}(\theta_1)$ whereas
Lemma \ref{sbeta} that
$\mathcal{N}(\beta)\subset\mathcal{N}(\theta_1)$.
Then \eqref{nucleogamma} gives
$$
\mathcal{N}(\gamma)=\mathcal{N}(\beta)\cap\mathcal{N}(\theta_1)
\cap\mathcal{N}(\theta_2)=\mathcal{N}(\beta)\cap\mathcal{N}(\theta_2).
$$
Hence, we have
$$
2n\geq\dim(\mathcal{N}(\beta)+\mathcal{N}(\theta_2))
=\nu(\beta)+\nu(\theta_2)
-\dim\mathcal{N}(\beta)\cap\mathcal{N}(\theta_2)
=\nu(\beta)+\nu(\theta_2)-\nu(\gamma)
$$
and since $\nu(\beta)=2n-2$, we conclude that 
$\nu(\gamma)\geq 2n-2p$.\qed

\begin{remark}\po{\em The estimate given by Lemma \ref{alglemma} 
is sharp. For instance, if we take as $\a$ in \eqref{gamma} and 
\eqref{beta} the second fundamental form of a product of real 
Kaehler hypersurfaces then the hypotheses of the lemma are satisfied 
and we have equality in the estimate.
}\end{remark}

\noindent{\em Proof of Theorem \ref{alglemma2}:} 
By Lemma \ref{alglemma} the vector subspace $\Sal(\gamma)$ 
is degenerate since, if otherwise, then by \eqref{T1} it is 
of the form $\Sal(\gamma)=W_1^{q,q}\subset W^{p,p}$ and then
Lemma \ref{alglemma} yields a contradiction.
Parts $(ii)$ and $(iv)$ of Proposition \ref{even}
give, respectively, that $\dim\Omega\geq 2$ and that 
$\gamma_\Omega=\pi_{\Omega\oplus\Omega}\circ\gamma$ satisfies
$\Sal(\gamma_\Omega)=\Sal(\gamma)\cap\Sal(\gamma)^\perp$.
Hence
$$
0=\lp\gamma_\Omega(X,Y),\gamma_\Omega(Z,T)\rp
=\<\a_\Omega(X,Y),\a_\Omega(Z,T)\>
-\<\a_\Omega(X,JY),\a_\Omega(Z,JT)\>
$$
for any $X,Y,Z,T\in V^{2n}$. Thus the complex structure 
$\J\in\text{End}(\Sal(\a_\Omega))$ defined by 
$\J\a_\Omega(X,Y)=\a_\Omega(X,JY)$ is an isometry. 
Part $(iv)$ of Proposition~\ref{even} gives that
$\Sal(\a_\Omega)=\Omega$ and hence $\J\in\text{Aut}(\Omega)$
is a complex structure. In particular, we have that 
$\a_\Omega$ is pluriharmonic. Then
$\Sal(\beta)\subset P\oplus P$ and hence $\beta=\beta_P
=\pi_{P\times P}\circ\beta$. 
Thus if $\gamma_P=\pi_{P\times P}\circ\gamma$ and
$\gamma=\gamma_\Omega\oplus\gamma_P$ then
$$
\lp\gamma_P(X,Y),\beta_P(Z,T)\rp
=\lp\gamma(X,Y),\beta(Z,T)\rp\;\;\mbox{for any}\;\;X,Y,Z,T\in V^{2n}.
$$
Hence $\gamma_P$ and $\beta_P$ 
satisfy the condition \eqref{productflat2}. Since $\gamma$ 
is flat and the bilinear form $\gamma_\Omega$ is null then 
also $\gamma_P$ is flat. Then by \eqref{T1} we have that
$\Sal(\gamma_P)=W^{q_1,q_1}$ and the remaining of the proof 
follows from Lemma \ref{alglemma}.\vspace{2ex}\qed

We now  prove the results stated in the introduction.
\vspace{2ex}

\noindent\emph{Proof of Theorem \ref{main}:} Let the bilinear forms
$\gamma,\beta\colon T_{x_0}M\times T_{x_0}M\to N_1^f(x_0)\oplus N_1^f(x_0)$ 
be defined  by \eqref{gamma} and \eqref{beta} in terms of the second 
fundamental form $\a$ of $f$ at $x_0\in M^{2n}$.  We endow
$N_1^f(x_0)\oplus N_1^f(x_0)$ with the inner product
defined by
$$
\lp(\xi,\bar\xi),(\eta,\bar\eta)\rp
=\<\xi,\eta\>_{N_1^f(x_0)}-\<\bar\xi,\bar\eta\>_{N_1^f(x_0)}.
$$
We claim that $\gamma$ and $\beta$ are flat and that 
\eqref{productflat2} holds. For a Kaehler manifold
it is a standard fact that the curvature tensor 
$x\in M^{2n}$ satisfies $R(X,Y)JZ=JR(X,Y)Z$ for any $X,Y,Z\in T_xM$. 
From this and the Gauss equation for $f$ we obtain
\begin{align*}
\lp&\gamma(X,T),\gamma(Z,Y)\rp= 
\<\a(X,T),\a(Z,Y)\>-\<\a(X,JT),\a(Z,JY)\>\\
&=\<R(X,Z)Y,T\>+\<\a(X,Y),\a(Z,T)\>-\<R(X,Z)JY,JT\>
-\<\a(X,JY),\a(Z,JT)\>\\
&=\lp\gamma(X,Y),\gamma(Z,T)\rp.
\end{align*}
Since $\gamma$ satisfies \eqref{T1} then using \eqref{mathcalT} 
we have
\begin{align*}
\lp\gamma(X,T),\beta(Z,Y)\rp
&=\lp\gamma(X,T),\gamma(Z,Y)\rp+\lp\gamma(X,T),\gamma(JZ,JY)\rp\\
&=\lp\gamma(X,Y),\gamma(Z,T)\rp+\lp\gamma(X,T),\T\gamma(JZ,Y)\rp\\
&=\lp\gamma(X,Y),\gamma(Z,T)\rp+\lp\gamma(X,JT),\gamma(JZ,Y)\rp\\
&=\lp\gamma(X,Y),\gamma(Z,T)\rp+\lp\gamma(X,Y),\gamma(JZ,JT)\rp\\
&=\lp\gamma(X,Y),\beta(Z,T)\rp.
\end{align*}
Using \eqref{productflat2} and since  $\beta(JX,Y)=-\beta(X,JY)$ 
from \eqref{beta} then 
\begin{align*}
&\lp\gamma(JX,JY),\beta(Z,T)\rp
=\lp\gamma(JX,T),\beta(Z,JY)\rp=
-\lp\gamma(JX,T),\beta(JZ,Y)\rp\\
&=-\lp\gamma(JX,Y),\beta(JZ,T)\rp=\lp\gamma(JX,Y),\beta(Z,JT)\rp
=\lp\gamma(JX,JT),\beta(Z,Y)\rp.
\end{align*}
Then by \eqref{beta} we have 
\begin{align*}
&\lp\beta(X,Y),\beta(Z,T)\rp
=\lp\gamma(X,Y),\beta(Z,T)\rp+\lp\gamma(JX,JY),\beta(Z,T)\rp\\
&=\lp\gamma(X,T),\beta(Z,Y)\rp+\lp\gamma(JX,JT),\beta(Z,Y)\rp
=\lp\beta(X,T),\beta(Z,Y)\rp
\end{align*}
and the claim has been proved. The proof now follows from 
Theorem \ref{alglemma2} since we have that
$\Delta_c(x_0)=\mathcal{N}(\gamma)$
and $Q(x_0)=\Omega$.\qed

\begin{remark}\label{ag}\po {\em  The proof of Theorem \ref{main} 
in \cite{CG} makes use of Theorem $1$ in \cite{DF1} but that result
does not hold for $p=6$ as was clarified in \cite{DF2}. Nevertheless, 
it was also established in \cite{DF2} that Theorem $1$ in \cite{DF1} 
still holds for $p=6$ under a slightly stronger assumption which 
happens to be satisfied in our case of real Kaehler submanifolds. 
}\end{remark}

\noindent\emph{Proof of Theorem \ref{cor}:} Theorem \ref{main}
and the Gauss equation for $f$ give
\begin{align*}
K(X,JX)
&=\<\a(X,X),\a(JX,JX)\>-\|\a(X,JX)\|^2\\
&=-\|\a_Q(X,X)\|^2-\|\a_Q(X,JX)\|^2\leq 0
\end{align*}
for any $X\in\mathcal{N}(\a_P)$.\qed
\vspace{3ex}

\noindent  Marcos Dajczer is partially 
supported by the grant PID2021-124157NB-I00 funded 
by MCIN/AEI/10.13039/501100011033/ ``ERDF A way of making Europe",
Spain, and is also supported by Comunidad Aut\'{o}noma de la Regi\'{o}n 
de Murcia, Spain, within the framework of the Regional Programme 
in Promotion of the Scientific and Technical Research (Action Plan 2022), 
by Fundaci\'{o}n S\'{e}neca, Regional Agency of Science and Technology, 
REF, 21899/PI/22.

\inputencoding{utf8}
{\renewcommand{\baselinestretch}{1}
\hspace*{-30ex}\begin{tabbing}
\indent \= Sergio J. Chion Aguirre \hspace{25ex} 
Marcos Dajczer\\
\>CENTRUM Católica Graduate Business School, \hspace{0.5ex} 
Departamento de Matemáticas\\
\> Lima, Perú, \hspace{36.9ex}
Universidad de Murcia, \\
\>Pontificia Universidad Catolica del Perú,  \hspace{8ex} 
Campus de Espinardo,\\\
\> Lima, Peru,\hspace{37.7ex} 
E-30100 Espinardo, Murcia, Spain\\
\>sjchiona@pucp.edu.pe \hspace{18.25ex}   \hspace{8ex} 
marcos@impa.br
\end{tabbing}}
\end{document}